\documentclass[10pt,a4paper]{article}
\usepackage[english]{babel}
\usepackage{lipsum}
\usepackage{commath}
\usepackage{graphicx}
\usepackage{enumitem}

\usepackage{float}
\usepackage{amsthm,amssymb,mathrsfs,amsmath}
\usepackage[a4paper,bindingoffset=0.2in,
left=0.5in,right=0.5in,
top=1in,bottom=1in,footskip=.25in]{geometry}
\title{Haar wavelets collocation on a class of Emden-Fowler equation via Newton's quasilinearization and Newton-Raphson techniques}
\author{Amit K. Verma$^a$, Narendra Kumar$^b$\\\small{\it{$^{a,b}$ Department of Mathematics,}} \\\small{\it{Indian Institute of Technology Patna,}}\\\small{\it{ Bihta, Patna 801103, (BR) India.}}}
\date{\today}
\theoremstyle{definition}

\newtheorem{theorem}{Theorem}[section]

\begin{document}

\maketitle
\begin{abstract}
In this paper we have considered generalized Emden-Fowler equation,
\begin{equation*}
y''(t)+\sigma t^\gamma y^\beta(t)=0, ~~~~~~~~t \in ]0,1[
\end{equation*}
subject to the following boundary conditions
\begin{equation*}
y(0)=1,~y(1)=0;~~\&~~y(0)=1,~y'(1)=y(1),
\end{equation*}  
where $\gamma,\beta$ and $\sigma$ are real numbers, $\gamma<-2$, $\beta>1$. We propsoed to solve the above singular nonlinear BVPs with the aid of Haar wavelet coupled with quasilinearization approach as well as Newton-Raphson approach. We have also considered the special case of Emden-Fowler equation ($\sigma=-1$,$\gamma=\frac{-1}{2}$ and $\beta=\frac{3}{2}$) which is popularly, known as Thomas-Fermi equation. We have analysed different cases of generalised Emden-Fowler equation and compared our results with existing results in literature. We observe that small perturbations in initial guesses does not affect the the final solution significantly. 
\end{abstract}
\section{Introduction}\label{intro}

The following singular boundary value problem,
\begin{equation*}
y''(t)= t^{-\frac{1}{2}} y^\frac{3}{2}(t), \quad\quad t \in ]0,1[,
\end{equation*}
with the following boundary conditions
\begin{equation*}
y(0)=1,~~y(1)=0,
\end{equation*}
is known as Thomas-Fermi equation (\cite{Numerical1998}). This equation was studied by Thomas \cite{Thomas1927} and Fermi \cite{Fermi1927}. One can impose three types of boundary conditions, which reflects three different states of atoms: \newline
a) ionized atom 
\begin{equation*}
y(0)=1,~~y(a)=0,
\end{equation*}
b) neutral atom
\begin{equation*}
y(0)=1,~~ay'(a)=y(a),
\end{equation*}
c) isolated neutral atom
\begin{equation*}
y(0)=1,~~\lim_{x\to\infty} y(x)=0.
\end{equation*}
Thomas-Fermi equation may be written in a generalised form known as Emden-Fowler equation (\cite{Numerical1998,Numerical1996})
\begin{eqnarray}
\label{EF}&&y''(t)+\sigma t^\gamma y^\beta(t)=0,\quad t \in ]0,1[,\\
\label{EFBC1}&&y(0)=1;~~y(1)=0,\\
\label{EFBC2}&&y(0)=1;~~y(1)=y'(1),
\end{eqnarray}
where $\gamma,\beta$ and $\sigma$ are real numbers, $\gamma<-2$, $\beta>1$. Mehta et al. \cite{MEHTA1971} studied it by using hypergeometric function and studied the Emden-Fowler equations. Wong \cite{JSWW1975} presented survey on Emden-Fowler equations. Rosenau \cite{ROSENAU1984} by simple change of variable generated a one parameter family of integrable Emden-Fowler equations. P M Lima \cite{Numerical1996, LIMA1999, LIMA2009} analyzed numerical methods based on finite differences/shooting methods for solving Emden-Fowler equation. Mohammadi et al. \cite{Mohammadi2019} used an iterative method to the singular and nonlinear fractional partial differential version of Emden-Fowler equations. Haar wavelets operational matrix of fractional integration is be applied to compute the solution with the Picard technique. 

Pikulin \cite{Pikulin2019} considered the Emden-Fowler equations on the half-line and on the interval and by assuming that the exponent in the coefficient of the nonlinear term is rational, he derived new parametric representations. 

In this paper we consider \eqref{EF} subject to the boundary conditions \eqref{EFBC1} and \eqref{EFBC2}.  We develop two approaches. In the first approach we use Newton's quasilinearization and arrive at an iteration scheme which is linear at each iteration. At each iteration we apply Haar wavelet collocation approach. After certain iterations the solutions does not vary, so we stop and the solution that we get at this stage is the solution of the Emden-Fowler equation. In the second approach we first use Haar wavelet collocation and arrive at the set of nonlinear equations. We then solve the system of nonlinear equation to get the solution of the Emden-Fowler equation. To the best of our knowledge the results do not exists in literature. 

The paper is divided into the following sections. In section \ref{intro} we discuss the literature survey. In section \ref{SecHaar} we discuss basics of Haar Wavelets. In section \ref{Method} we discuss about proposed methodologies. Section \ref{Numerics} is devoted to computations. Finally in section \ref{Conclusion} we conclude the paper.

\section{Haar Wavelet}\label{SecHaar}

\cite[pp.7-10]{Haar} Let us consider the interval $0\leq x \leq 1$. Let us define $M=2^J$,where $J$ is maximum level of resolution. Let us divide $[0,1]$ into $2M$ subintervals of equal length $\Delta x=\frac{1}{2M}$. The wavelet number $i$ is calculated by $i=m+k+1$,here $j=0,1,\cdots,J$ and $k=0,1,\cdots,m-1$(here $m=2^j$). The Haar wavelet's mother wavelet function is defined as,
\begin{equation} \label{eq1}
h_i(x)=
\begin{cases} 
     1,  & {\alpha_1(i)\leq x< \alpha_2(i)}, \\
    -1,  & {\alpha_2(i)\leq x< \alpha_3(i)}, \\
      0, & \quad \text{else},
   \end{cases}
\end{equation}
where,
\begin{equation} \label{eq2} \quad \alpha_1(i)=2k \mu \Delta x , \quad \alpha_2(i)=(2k+1) \mu \Delta x , \quad \alpha_3(i)=2(k+1) \mu \Delta x , \quad \mu =\frac{M}{m},
\end{equation} \newline
if $i>2$. For $i=1$, it is defined as,
\begin{equation} \label{eq3}
h_1(x)=
\begin{cases} 
     1,  & {0 \leq x \leq 1}, \\
      0, & \quad \text{else}.
   \end{cases}
\end{equation}
For $i=2$, 
\begin{equation} \label{eq4}
\quad \alpha_1(2)=0 , \quad \alpha_2(2)=0.5 , \quad \alpha_3(2)=1.
\end{equation}
The width of the $i^{th}$ wavelet is,
\begin{center}
{$\alpha_3(i)-\alpha_1(i)= 2 \mu \Delta x = 2^{-j}$}.
\end{center}
The integral $P_{v,i}(x)$ is defined as,
\begin{equation} \label{eq5}
P_{v,i}(x)=\int^x_0 \int^x_0 \cdots \int^x_0 h_i(t) dt^v = \frac{1}{(v-1)!} \int^x_0 (x-t)^{v-1} h_i(t) dt,
\end{equation}
here $v$ is the order of integration. \newline
Using (\ref{eq1}), we will calculate these integrals analytically and by doing it we obtain,
\begin{equation} \label{eq6}
P_{v,i}(x)=
\begin{cases} 
     0,  & {x< \alpha_1(i)}, \\
    \frac{1}{v!}[x-\alpha_1(i)]^v,  & {\alpha_1(i)\leq x \leq \alpha_2(i)}, \\
    \frac{1}{v!}{[x-\alpha_1(i)]^v - 2[x-\alpha_2(i)]^v},  & {\alpha_2(i)\leq x \leq \alpha_3(i)},\\
      \frac{1}{v!}{[x-\alpha_1(i)]^v - 2[x-\alpha_2(i)]^v + [x-\alpha_3(i)]^v},  & { x > \alpha_3(i)},
   \end{cases}
\end{equation} 
for $i>1$. For $i=1$ we have $\alpha_1=0,\alpha_2=\alpha_3=1$ and,
\begin{equation} \label{eq7}
P_{v,1}(x)=\frac{x^v}{v!}.
\end{equation}
For computation by using Haar wavelet, we make use of the method of collocation. Here, the grid points are defined by,
\begin{equation} \label{eq8}
\tilde x_c=c \Delta x,\quad c=0,1,\cdots,2M,
\end{equation}
and collocation points are defined by,
\begin{equation} \label{eq9}
x_c=0.5(\tilde x_{c-1}+\tilde{x_c}),\quad c=1,\cdots,2M.
\end{equation}
We substitute $x \rightarrow x_c$ in (\ref{eq1}),(\ref{eq6}),(\ref{eq7}).

For computational point of view, we introduce the Haar matrices $H,P_1,P_2,\cdots ,P_v$. The order of these matrices are $2M \times 2M$. These matrices are defined by,
\begin{center}
$H(i,c)=h_i(x_c), P_v(i,c)=p_{v,i}(x_c) ~~~v=1,2,\cdots.$
\end{center} 
For $J=1$, the matrices $H$, $P_1$  and $P_2$ are defined by,
\begin{center}
$H=
\begin{bmatrix}
1&1&1&1 \\ 1&1&-1&-1 \\ 1&-1&0&0 \\ 0&0&1&-1
\end{bmatrix}$,~~ $P_1= \frac{1}{8} \begin{bmatrix}
1&3&5&7 \\ 1&3&3&1 \\ 1&1&0&0 \\ 0&0&1&1
\end{bmatrix}$,~~ $P_2= \frac{1}{128} \begin{bmatrix}
1&9&25&49 \\ 1&9&23&31 \\ 1&7&8&8 \\ 0&0&1&7
\end{bmatrix}.$ 
\end{center}

\section{Haar Wavelet Collocation Method}\label{Method}

In this section we develop two methods for solving generalized Emden-Fowler equation.

\subsection{Quasilinearization Approach}

By quasilinearization (\cite{An2013, Higher2019}) generalized Emden-Fowler equation is linearized then at each iteration Haar wavelet collocation method is used to compute the solution. After some iteration the solution converges to the solution of the Emden-Fowler equation. Consider the problem considerd by Lima \cite{Numerical1996},

\begin{equation} \label{eq10}
y''(t)+\sigma t^\gamma y^\beta(t)=0, ~~~~~~~~t \in ]0,1[
\end{equation}
\begin{equation} \label{eq11}
y(0)=1;~~y(1)=0,
\end{equation}
where $\gamma,\beta$ and $\sigma$ are real numbers, $\gamma<-2$, $\beta>1$. Now applying quasilinearization in above equation,

\begin{equation} \label{eq12}
y_{r+1}''(t)=\sigma t^\gamma [-y_r^\beta(t) + (y_{r+1}(t)-y_r(t))(-\beta y_r^{\beta -1}(t))],
\end{equation}
\begin{equation} \label{eq13}
y_{r+1}(0)=1; ~~~y_r(1)=0.
\end{equation}
To apply Haar wavelet method \cite{Haar_solving_HODE2008}, we assume that

\begin{equation} \label{eq14}
y_{r+1}''(t)=\sum_{i=1}^{2M} a_i h_i(t),
\end{equation}
where $a_i$ are wavelet coefficients. Now integrate above equation from $0$ to $t$ twice we get,

\begin{equation} \label{eq15}
y_{r+1}'(t)=\sum_{i=1}^{2M} a_i P_{1,i}(t) + y_{r+1}'(0),
\end{equation}
\begin{equation} \label{eq16}
y_{r+1}(t)=\sum_{i=1}^{2M} a_i P_{2,i}(t) + ty_{r+1}'(0)+ y_{r+1}(0).
\end{equation}

Now apply boundary conditions (\ref{eq13}) in (\ref{eq16}) and find $y_{r+1}'(0)$ then substitute (\ref{eq16}),(\ref{eq14}) in (\ref{eq12}) and we get
\begin{equation} \label{eq17}
\sum_{i=1}^{2M} a_i [h_i(t)+\sigma \beta t^{\gamma}y_r^{\beta-1}(t)(P_{2,i}(t)-tP_{2,i}(1))]=\sigma (\beta-1) t^{\gamma}y_r^{\beta}(t)-\sigma \beta (1-t) t^{\gamma} y_r^{\beta-1}(t),
\end{equation}

Now, we can solve above equation by collocation method.

\subsection{Newton Raphson Method}

\cite[pp.33]{Haar} Consider a nonlinear ordinary differential equation,

\begin{equation} \label{eq18}
\psi (x,y,y',y'')=0, ~~~~x\in[0,1],
\end{equation}
where $\psi$ is a nonlinear function. Now again applying Haar wavelet method \cite{Haar_solving_HODE2008}, let us assume
\begin{equation} \label{eq19}
y''(t)=\sum_{i=1}^{2M} a_i h_i(t),
\end{equation}
and integrate above equation $2$-times. Obtained results are discretized by collocation method then substitute these results in (\ref{eq18}). By doing this we will get a system of non-linear equations,

\begin{equation} \label{eq20}
\phi_c (a_1,a_2,\cdots,a_{2M})=0, ~~~c=1,2,\cdots,2M.
\end{equation} 
To solve above equation here we are using Newton's method and using this method wavelet coefficients ($a_i$) can be calculated.

\subsection{Convergence} The convergence of the proposed methods is an immediate consequence of the results stated by Majak et al. \cite{BS2015}.
\begin{theorem} (\cite{BS2015})
Let us assume that $g(x)=\frac{d^2 y}{dx^2}\in L^2(\mathbb{R})$ is a continuous function on $[0,1]$ and its first derivative is bounded for all $x\in [0,1]$, there exists $\xi$ such that $\left|\frac{dg}{dt}\right|\leq \xi$, then both approaches discussed in this paper converge.
\end{theorem}
Proof. The proof is similar to the proof given in \cite{Higher2019}.
\section{Numerical Illustration} \label{Numerics}

\subsection{Example 1}

Consider (\ref{eq10}), (\ref{eq11}) with $\sigma=-1$,$\gamma=\frac{-1}{2}$ and $\beta=\frac{3}{2}$. Solving by quasilinearization method (subsection $3.1$) substitute $\sigma$,$\gamma$ and $\beta$ and replace $t \rightarrow t_c$ $\forall$ $c=1,2,\cdots,2M$ in (\ref{eq17}) where $t_c$ are collocation points (section $2$). Now we get (\ref{eq17}) in terms of $H$,$P_2$,$A$ and $Y_r$ where $A=[a_1,a_2,\cdots,a_{2M}]$ and $Y_r=[y_r(t_1),y_r(t_2),\cdots,y_r(t_{2M})]$ are row vectors from which unknowns $[a_1,a_2,\cdots,a_2M]$ are calculated.

\begin{figure}[H]
\begin{center}
\includegraphics[scale=0.3]{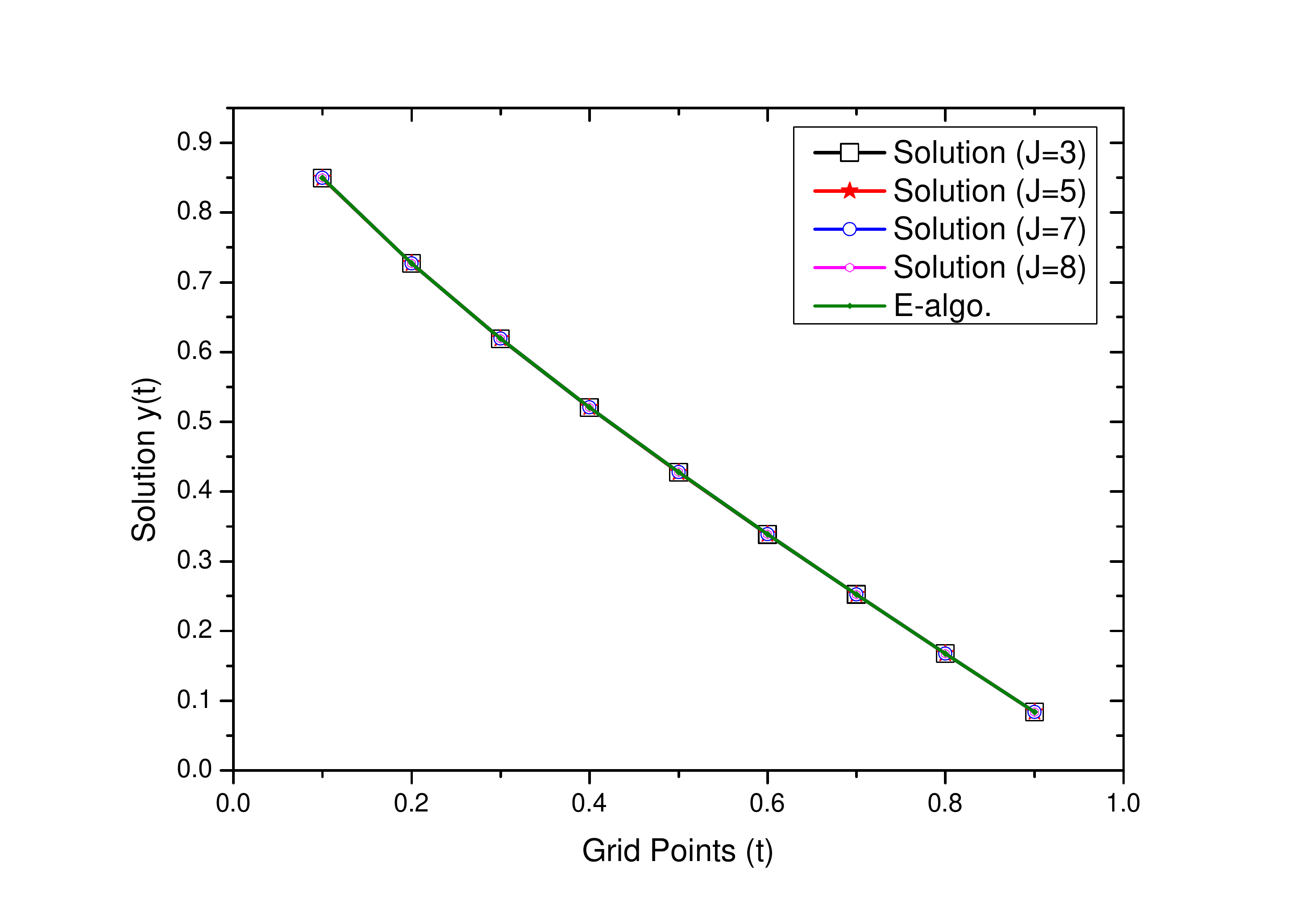}
\caption{Graph of Quasilinearization method for $J=3$, $J=5$, $J=7$, $J=8$ and E-algo  (\cite{Numerical1996}).}\label{ppr1fig1}
\end{center}
\end{figure}

Taking initial guess $[y_0(t_1),y_0(t_2),\cdots,y_0(t_{2M})]$ and by Haar wavelet collocation method, we will get required solution. Taking initial guess $[0,0,\cdots,0]$ solution is given in Table \ref{PPR1Table1} and graph for $J=3$, $J=5$, $J=7$ and $J=8$ is given in figure \ref{ppr1fig1}. Taking different initial guess $[0.01,0.01,\cdots,0.01]$ and $[0.1,0.1,\cdots,0.1]$ solution still remain same.

Now by Newton's method (subsection $3.2$), assume 

\begin{equation} \label{eq21}
y''(t)=\sum_{i=1}^{2M} a_i h_i(t),
\end{equation}
where $a_i$ are wavelet coefficients. Now integrate above equation from $0$ to $t$ twice we get,

\begin{equation} \label{eq22}
y'(t)=\sum_{i=1}^{2M} a_i P_{1,i}(t) + y'(0),
\end{equation}
\begin{equation} \label{eq23}
y(t)=\sum_{i=1}^{2M} a_i P_{2,i}(t) + ty'(0)+ y(0).
\end{equation}
Now apply boundary conditions (\ref{eq13}) in (\ref{eq23}) and we get $y(t)$ as,
\begin{equation} \label{eq24}
y(t)=1-t+\sum_{i=1}^{2M} a_i (P_{2,i}(t) - tP_{2,i}(1)).
\end{equation} 
Substituting (\ref{eq24}),(\ref{eq21}) in (\ref{eq18}), we get system of non-linear equations. Thus we arrive at (\ref{eq20}).  To solve the non-linear equation we use Newton's method to calculate wavelet coefficients ($a_i$).
 
\begin{figure}[H]
\begin{center}
\includegraphics[scale=0.3]{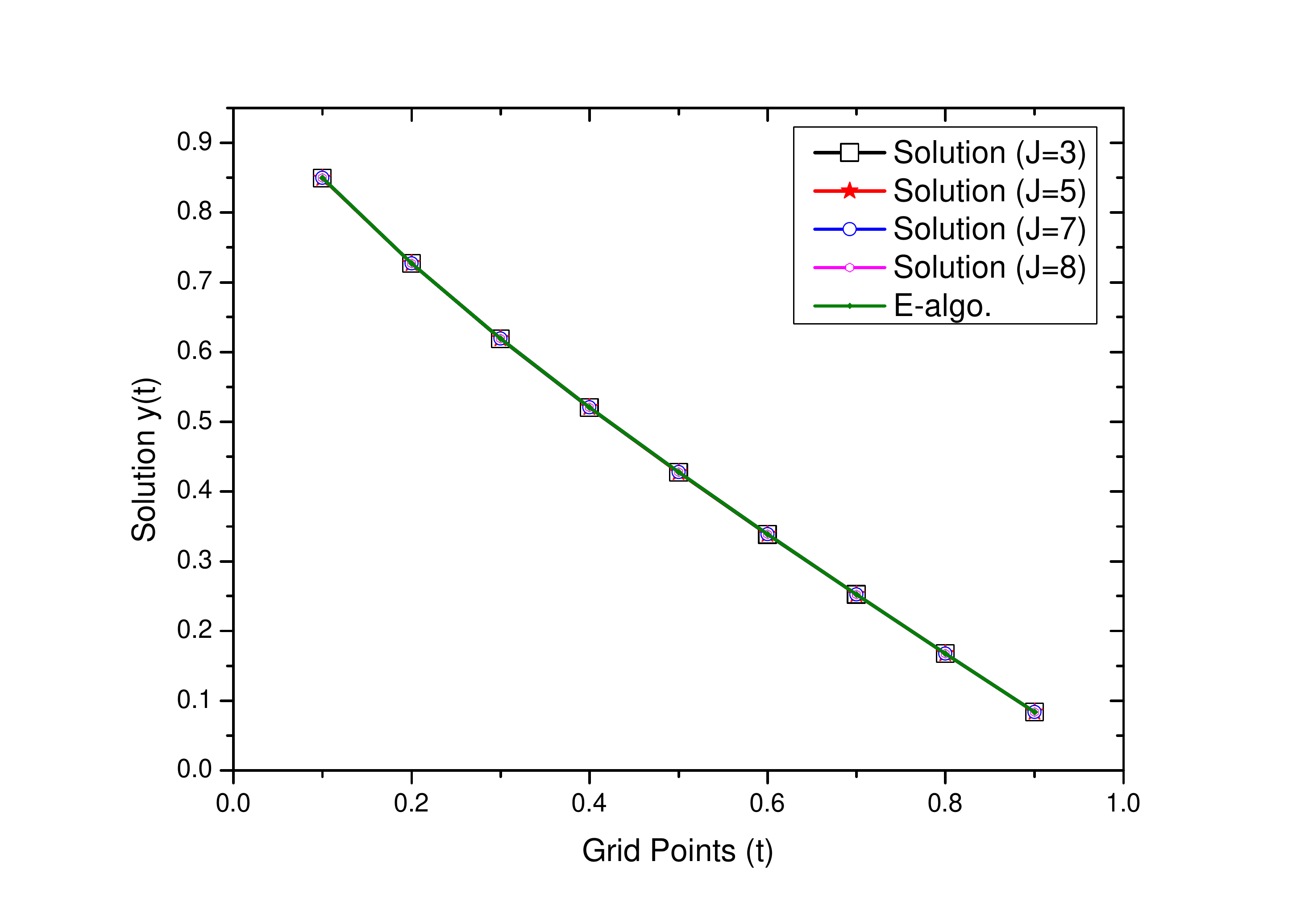}
\caption{Graph of Newton's method for $J=3$, $J=5$, $J=7$, $J=8$ and E-algo (\cite{Numerical1996}).}\label{ppr1fig2}
\end{center}
\end{figure}

For initial guess $[1,1,\cdots,1]$ computed solution is given in Table \ref{PPR1Table2} and graph for $J=3$, $J=5$, $J=7$ and $J=8$ is given in figure \ref{ppr1fig2}. Taking varying initial guesses to  $[1.01,1.01,\cdots,1.01]$ and $[1.1,1.1,\cdots,1.1]$ we observe that solution does not change significantly.

\begin{table}[H]											 
\centering											
\begin{center}											
\resizebox{10cm}{2.5cm}{											
\begin{tabular}	{|c | c|  c|  c|  c|  c|}		
\hline
$t$	&	$J=3,  r=3$	&	$J=5,  r=3$	&	$J=7,  r=3$	&	$J=8,  r=3$	&	E-algo.\cite[pp.262]{Numerical1996}	\\\hline
0.1	&	0.849094	&	0.849401	&	0.849464	&	0.84947	&	0.849474	\\\hline
0.2	&	0.726803	&	0.727162	&	0.727222	&	0.727228	&	0.727231	\\\hline
0.3	&	0.618867	&	0.619233	&	0.619286	&	0.619292	&	0.619294	\\\hline
0.4	&	0.520041	&	0.520361	&	0.520408	&	0.520412	&	0.520414	\\\hline
0.5	&	0.427227	&	0.427506	&	0.427544	&	0.427548	&	0.42755	\\\hline
0.6	&	0.338421	&	0.338651	&	0.338682	&	0.338685	&	0.338686	\\\hline
0.7	&	0.252197	&	0.252371	&	0.252395	&	0.252397	&	0.252398	\\\hline
0.8	&	0.16751	&	0.167631	&	0.167647	&	0.167648	&	0.167649	\\\hline
0.9	&	0.0836165	&	0.0836776	&	0.0836856	&	0.0836864	&	0.083686	\\\hline

\end{tabular}}								
\end{center}
\caption{\small{Solutions for $\sigma=-1$,$\gamma=\frac{-1}{2}$ and $\beta=\frac{3}{2}$ for $J=3$, $J=5$, $J=7$ and $J=8$ by Quasilinearization method}}	
\label{PPR1Table1}											
\end{table}

\begin{table}[H]											 
\centering											
\begin{center}											
\resizebox{10cm}{2.5cm}{										
\begin{tabular}	{|c | c|  c|  c|  c|  c|}		
\hline
$t$	&	$J=3,  r=3$	&	$J=5,  r=3$	&	$J=7,  r=3$	&	$J=8,  r=3$	&	E-algo.\cite[pp.262]{Numerical1996}	\\\hline
0.1	&	0.849094	&	0.849401	&	0.849464	&	0.84947	&	0.849474	\\\hline
0.2	&	0.726803	&	0.727162	&	0.727222	&	0.727228	&	0.727231	\\\hline
0.3	&	0.618867	&	0.619233	&	0.619286	&	0.619292	&	0.619294	\\\hline
0.4	&	0.520041	&	0.520361	&	0.520408	&	0.520412	&	0.520414	\\\hline
0.5	&	0.427227	&	0.427506	&	0.427544	&	0.427548	&	0.42755	\\\hline
0.6	&	0.338421	&	0.338651	&	0.338682	&	0.338685	&	0.338686	\\\hline
0.7	&	0.252197	&	0.252371	&	0.252395	&	0.252397	&	0.252398	\\\hline
0.8	&	0.16751	&	0.167631	&	0.167647	&	0.167648	&	0.167649	\\\hline
0.9	&	0.0836165	&	0.0836776	&	0.0836856	&	0.0836864	&	0.083686	\\\hline

\end{tabular}}								
\end{center}
\caption{\small{Solutions for $\sigma=-1$,$\gamma=\frac{-1}{2}$ and $\beta=\frac{3}{2}$ for $J=3$, $J=5$, $J=7$ and $J=8$ by Newton's method}}	
\label{PPR1Table2}											
\end{table}

\subsection{Example 2}

\begin{figure}[H]
\begin{center}
\includegraphics[height=7 cm,width=9 cm]{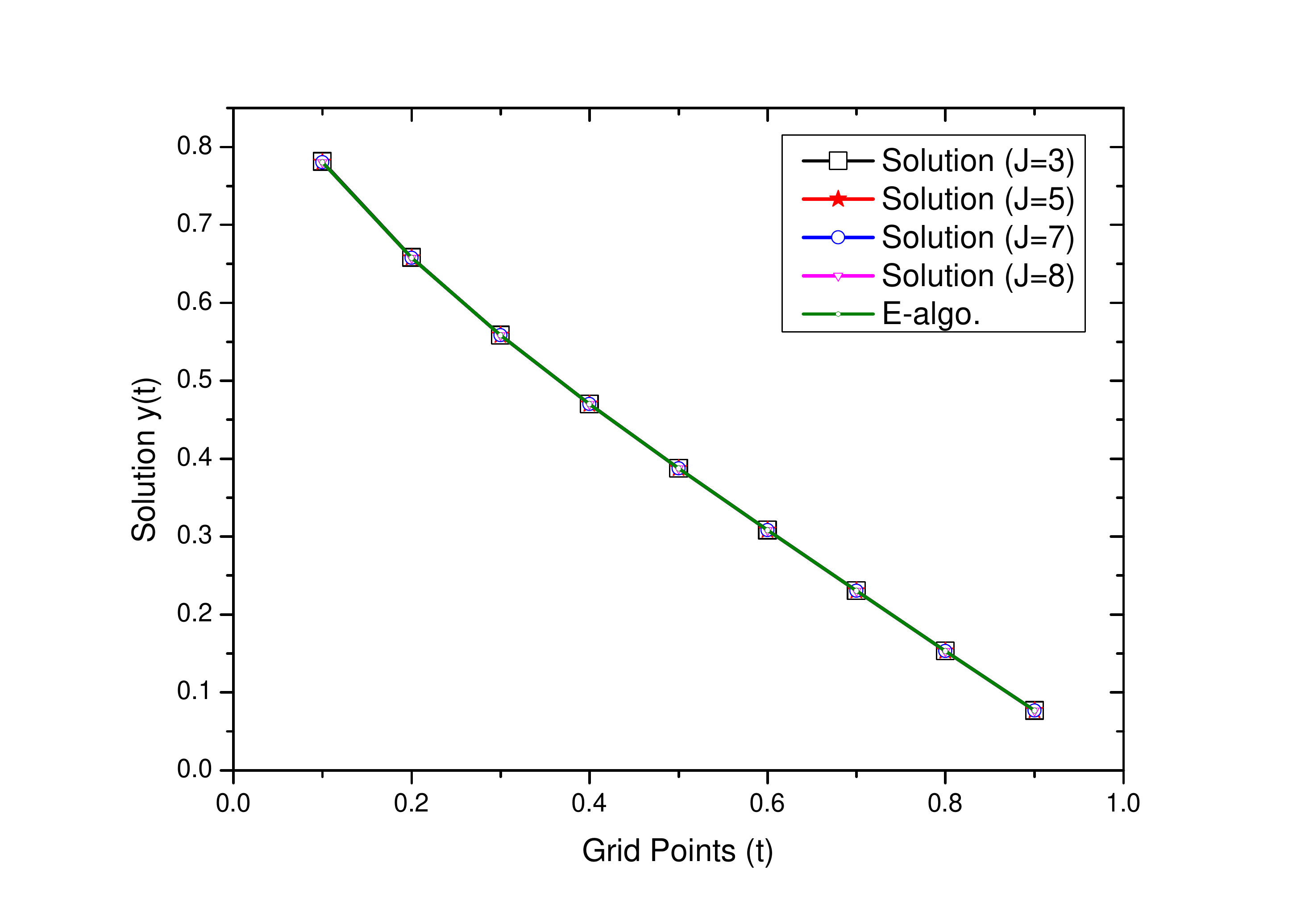}
\caption{Graph of Quasilinearization method for $J=3$, $J=5$, $J=7$, $J=8$ and E-algo  (\cite{Numerical1996}).}\label{ppr1fig3}
\end{center}
\end{figure}

Consider (\ref{eq10}), (\ref{eq11}) with $\sigma=-1$,$\gamma=-1$ and $\beta=2$. Here we follow same procedure similar to example $1$. For Quasilinearization method we take $[0,0,\cdots,0]$ as initial guess. Computed solution is given in Table \ref{PPR1Table3} and graph for $J=3$, $J=5$, $J=7$ and $J=8$ is given in figure \ref{ppr1fig3}. Taking different initial guesses $[0.01,0.01,\cdots,0.01]$ and $[0.1,0.1,\cdots,0.1]$ we observe that solution does not vary which proves that method is stable.

\begin{figure}[H]
\begin{center}
\includegraphics[scale=0.3]{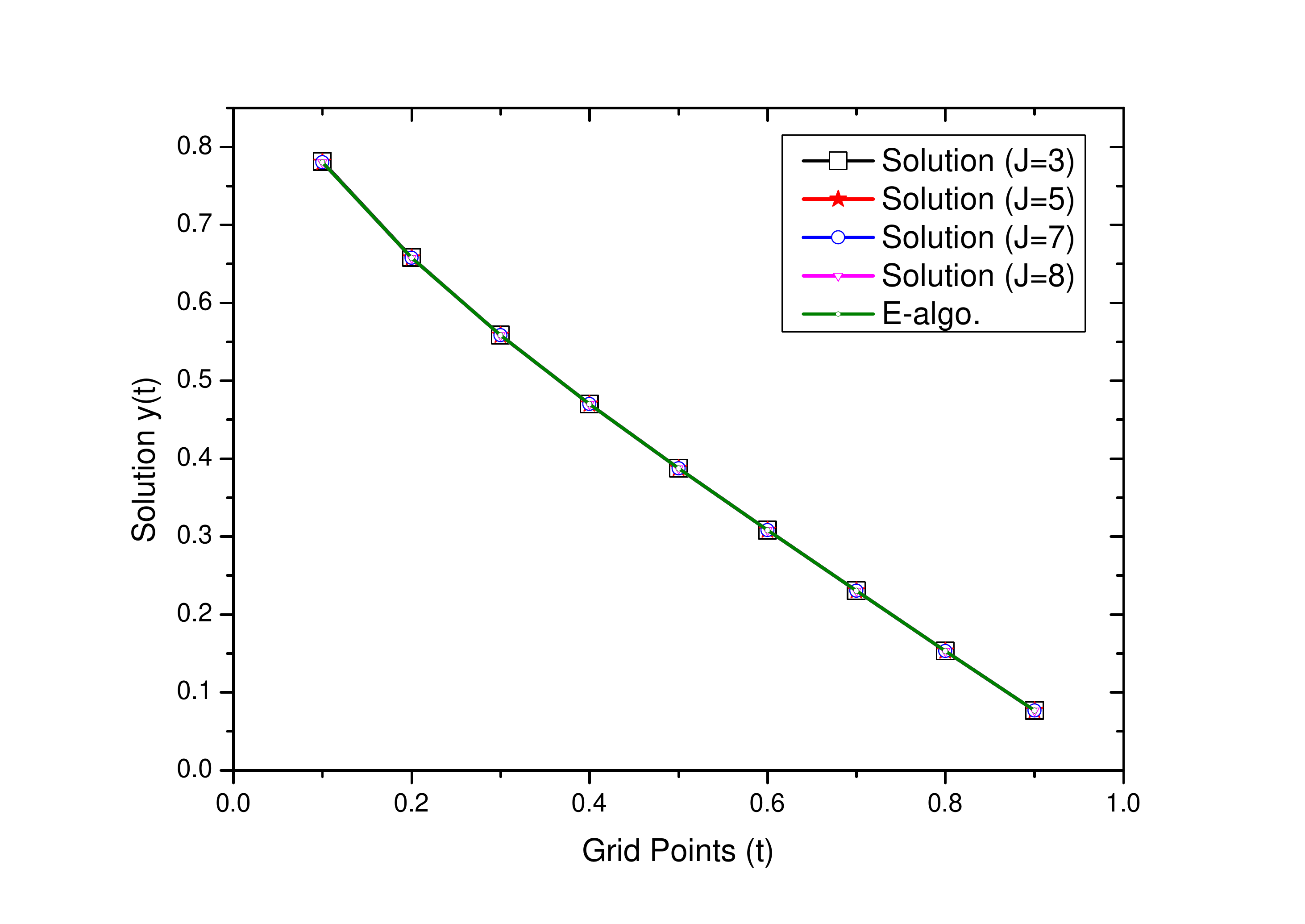}
\caption{Graph of Newton's method for $J=3$, $J=5$, $J=7$, $J=8$ and E-algo  (\cite{Numerical1996}).}\label{ppr1fig4}
\end{center}
\end{figure}

For Newton's method we take $[1,1,\cdots,1]$ as initial guess. solution is given in Table \ref{PPR1Table4} and graph for $J=3$, $J=5$, $J=7$ and $J=8$ is given in figure \ref{ppr1fig4}. Taking different initial guess $[1.01,1.01,\cdots,1.01]$ and $[1.1,1.1,\cdots,1.1]$ solution still remains same.

\subsection{Example 3}

\begin{figure}[H]
\begin{center}
\includegraphics[scale=0.3]{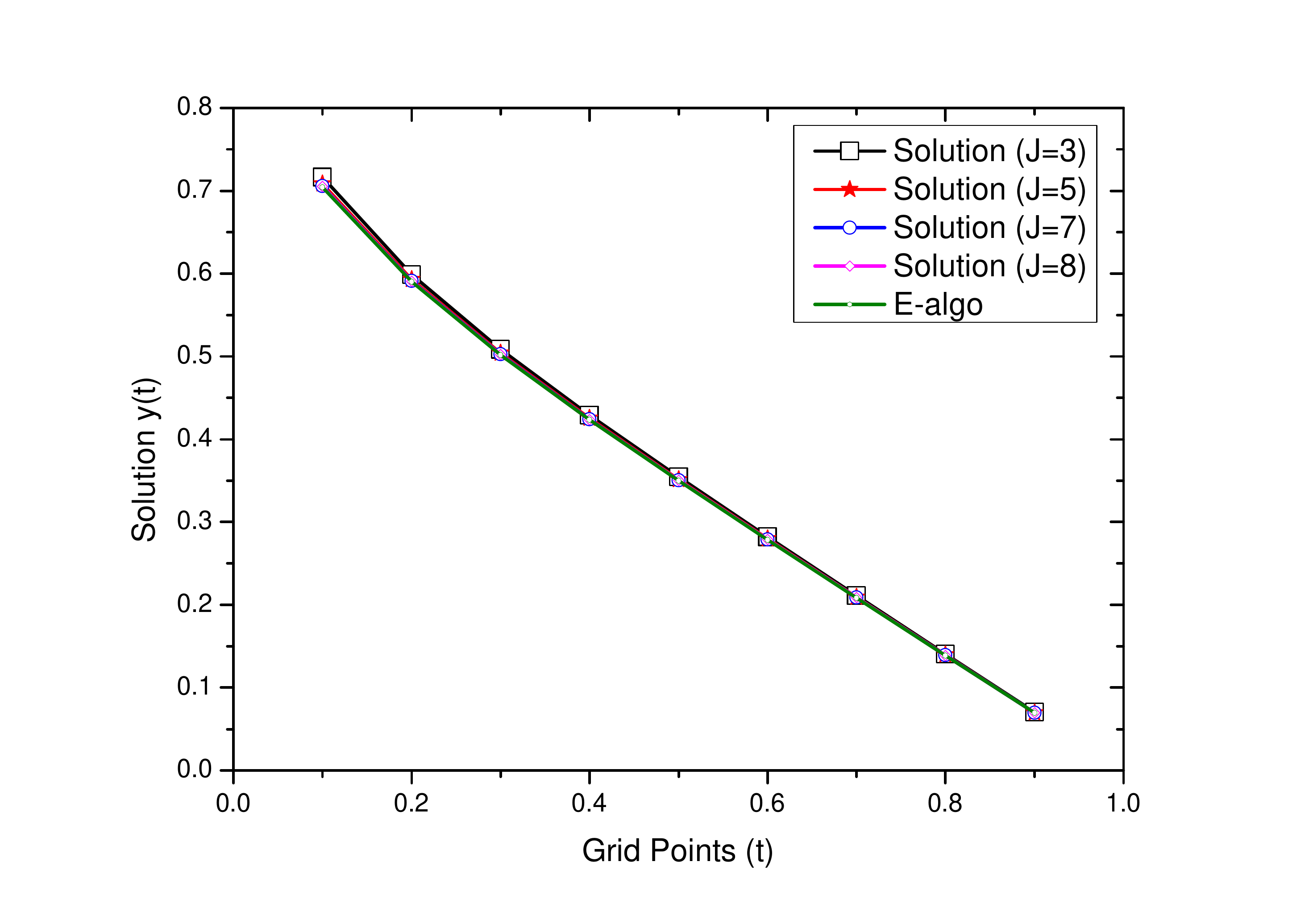}
\caption{Graph of Quasilinearization method for $J=3$, $J=5$, $J=7$, $J=8$ and E-algo  (\cite{Numerical1996}).}\label{ppr1fig5}
\end{center}
\end{figure}

Consider (\ref{eq10}), (\ref{eq11}) with $\sigma=-1$,$\gamma=\frac{-5}{4}$ and $\beta=\frac{9}{4}$. Here we follow same procedure similar to example $1$. For Quasilinearization method we take $[0,0,\cdots,0]$ as initial guess. Computed solution is displayed in Table \ref{PPR1Table5} and graph for $J=3$, $J=5$, $J=7$ and $J=8$ is given in figure \ref{ppr1fig5}. Small perturbations in initial guesses $[0.01,0.01,\cdots,0.01]$ and $[0.1,0.1,\cdots,0.1]$ does not change the solution significantly.

\begin{figure}[H]
\begin{center}
\includegraphics[scale=0.3]{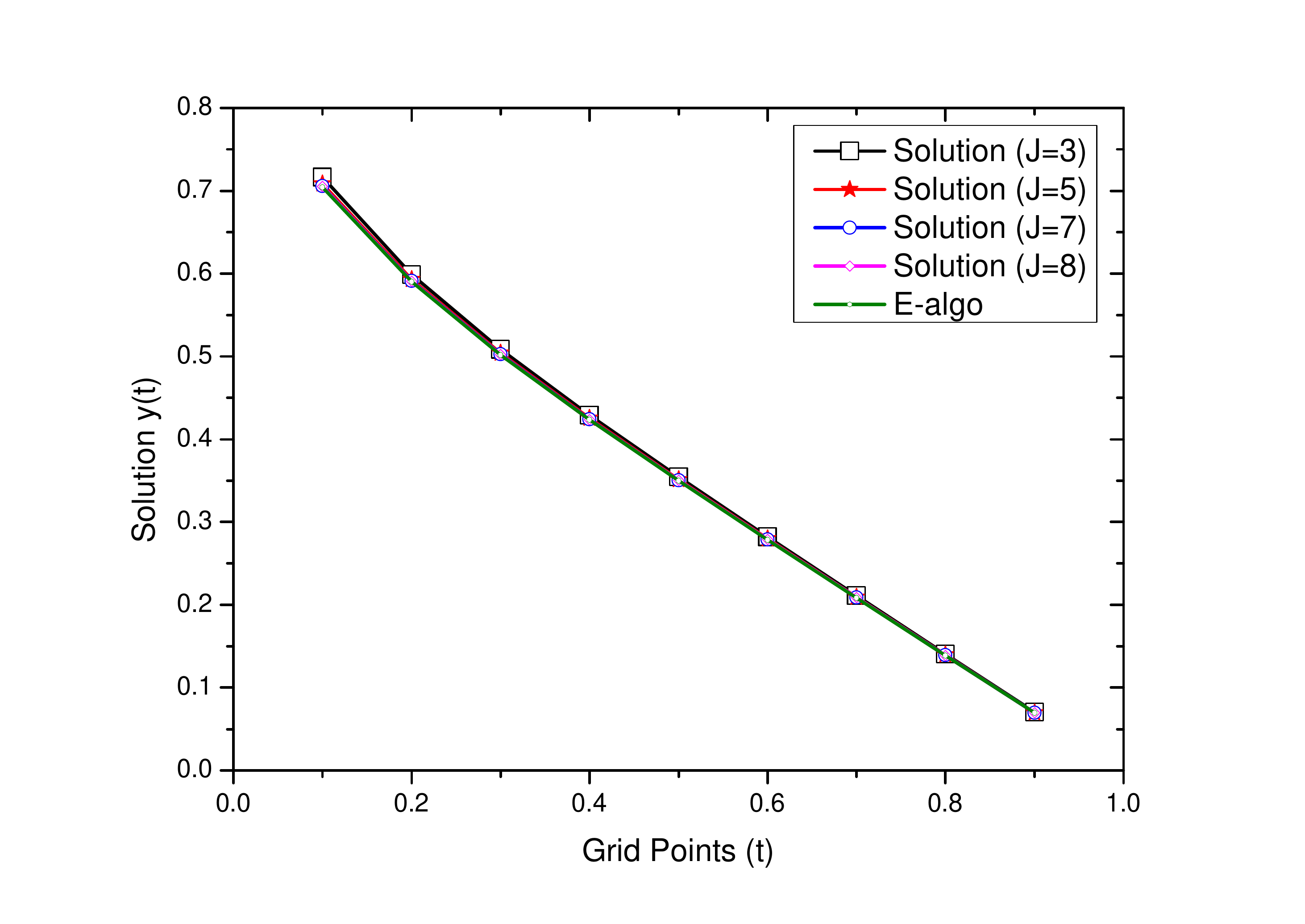}
\caption{Graph of Newton's method for $J=3$, $J=5$, $J=7$, $J=8$ and E-algo  (\cite{Numerical1996}).}\label{ppr1fig6}
\end{center}
\end{figure}

For Newton's method we take $[1,1,\cdots,1]$ as initial guess. Computed solution is given in Table \ref{PPR1Table6} and graph for $J=3$, $J=5$, $J=7$ and $J=8$ is given in figure \ref{ppr1fig6}. Taking different initial guess $[1.01,1.01,\cdots,1.01]$ and $[1.1,1.1,\cdots,1.1]$ we observe that solution remains same.

\begin{table}[H]											 
\centering											
\begin{center}											
\resizebox{10cm}{2.5cm}{										
\begin{tabular}	{|c | c|  c|  c|  c|  c|}		
\hline
$t$	&	$J=3, r=3$	&	$J=5, r=3$	&	$J=7, r=3$	&	$J=8, r=3$	&	E-algo.\cite[pp.264]{Numerical1996}	\\\hline
0.1	&	0.781336	&	0.780201	&	0.78013	&	0.780126	&	0.780125	\\\hline
0.2	&	0.657973	&	0.657498	&	0.657471	&	0.657469	&	0.657468	\\\hline
0.3	&	0.558596	&	0.558365	&	0.55835	&	0.558349	&	0.558348	\\\hline
0.4	&	0.470276	&	0.470119	&	0.470109	&	0.470109	&	0.470108	\\\hline
0.5	&	0.387687	&	0.387587	&	0.387581	&	0.387581	&	0.38758	\\\hline
0.6	&	0.308217	&	0.30815	&	0.308146	&	0.308146	&	0.308145	\\\hline
0.7	&	0.23039	&	0.230345	&	0.230342	&	0.230342	&	0.230342	\\\hline
0.8	&	0.153353	&	0.153328	&	0.153327	&	0.153326	&	0.153326	\\\hline
0.9	&	0.0766365	&	0.0766246	&	0.0766239	&	0.0766238	&	0.076623	\\\hline
\end{tabular}}								
\end{center}
\caption{\small{Solutions for $\sigma=-1$,$\gamma=-1$ and $\beta=2$ for $J=3$, $J=5$, $J=7$ and $J=8$ by Quasilinearization method}}	
\label{PPR1Table3}											
\end{table}

\begin{table}[H]											 
\centering											
\begin{center}											
\resizebox{10cm}{2.5cm}{
\begin{tabular}	{|c | c|  c|  c|  c|  c|}		
\hline
$t$	&	$J=3, r=3$	&	$J=5, r=3$	&	$J=7, r=3$	&	$J=8, r=3$	&	E-algo.\cite[pp.264]{Numerical1996}	\\\hline
0.1	&	0.781336	&	0.780201	&	0.78013	&	0.780126	&	0.780125	\\\hline
0.2	&	0.657973	&	0.657498	&	0.657471	&	0.657469	&	0.657468	\\\hline
0.3	&	0.558596	&	0.558365	&	0.55835	&	0.558349	&	0.558348	\\\hline
0.4	&	0.470276	&	0.470119	&	0.470109	&	0.470109	&	0.470108	\\\hline
0.5	&	0.387687	&	0.387587	&	0.387581	&	0.387581	&	0.38758	\\\hline
0.6	&	0.308217	&	0.30815	&	0.308146	&	0.308146	&	0.308145	\\\hline
0.7	&	0.23039	&	0.230345	&	0.230342	&	0.230342	&	0.230342	\\\hline
0.8	&	0.153353	&	0.153328	&	0.153327	&	0.153326	&	0.153326	\\\hline
0.9	&	0.0766365	&	0.0766246	&	0.0766239	&	0.0766238	&	0.076623	\\\hline

\end{tabular}}								
\end{center}
\caption{\small{Solutions for $\sigma=-1$,$\gamma=-1$ and $\beta=2$ for $J=3$, $J=5$, $J=7$ and $J=8$ by Newton's method}}	
\label{PPR1Table4}											
\end{table}

\begin{table}[H]											 
\centering											
\begin{center}											
\resizebox{10cm}{2.5cm}{										
\begin{tabular}	{|c | c|  c|  c|  c|  c|}		
\hline
$t$	&	$J=3, r=3$	&	$J=5, r=3$	&	$J=7, r=3$	&	$J=8, r=3$	&	E-algo.\cite[pp.265]{Numerical1996}	\\\hline
0.1	&	0.716249	&	0.708156	&	0.7057	&	0.70517	&	0.704396	\\\hline
0.2	&	0.598674	&	0.593057	&	0.591181	&	0.590768	&	0.590163	\\\hline
0.3	&	0.508462	&	0.504055	&	0.502524	&	0.502185	&	0.501688	\\\hline
0.4	&	0.428945	&	0.425335	&	0.424071	&	0.42379	&	0.423379	\\\hline
0.5	&	0.354347	&	0.351425	&	0.350393	&	0.350164	&	0.349827	\\\hline
0.6	&	0.282173	&	0.27987	&	0.279053	&	0.278872	&	0.278605	\\\hline
0.7	&	0.211155	&	0.20944	&	0.208831	&	0.208696	&	0.208497	\\\hline
0.8	&	0.140637	&	0.1395	&	0.139095	&	0.139005	&	0.138872	\\\hline
0.9	&	0.0702995	&	0.0697316	&	0.0695291	&	0.0694841	&	0.069418	\\\hline

\end{tabular}}								
\end{center}
\caption{\small{Solutions for $\sigma=-1$,$\gamma=\frac{-5}{4}$ and $\beta=\frac{9}{4}$ for $J=3$, $J=5$, $J=7$ and $J=8$ by Quasilinearization method}}	
\label{PPR1Table5}											
\end{table}

\begin{table}[H]											 
\centering											
\begin{center}											
\resizebox{10cm}{2.5cm}{											
\begin{tabular}	{|c | c|  c|  c|  c|  c|}		
\hline
$t$	&	$J=3, r=3$	&	$J=5, r=3$	&	$J=7, r=3$	&	$J=8, r=3$	&	E-algo.\cite[pp.265]{Numerical1996}	\\\hline
0.1	&	0.716249	&	0.708156	&	0.7057	&	0.70517	&	0.704396	\\\hline
0.2	&	0.598674	&	0.593057	&	0.591181	&	0.590768	&	0.590163	\\\hline
0.3	&	0.508462	&	0.504055	&	0.502524	&	0.502185	&	0.501688	\\\hline
0.4	&	0.428945	&	0.425335	&	0.424071	&	0.42379	&	0.423379	\\\hline
0.5	&	0.354347	&	0.351425	&	0.350393	&	0.350164	&	0.349827	\\\hline
0.6	&	0.282173	&	0.27987	&	0.279053	&	0.278872	&	0.278605	\\\hline
0.7	&	0.211155	&	0.20944	&	0.208831	&	0.208696	&	0.208497	\\\hline
0.8	&	0.140637	&	0.1395	&	0.139095	&	0.139005	&	0.138872	\\\hline
0.9	&	0.0702995	&	0.0697316	&	0.0695291	&	0.0694841	&	0.069418	\\\hline

\end{tabular}}								
\end{center}
\caption{\small{Solutions for $\sigma=-1$,$\gamma=\frac{-5}{4}$ and $\beta=\frac{9}{4}$ for $J=3$, $J=5$, $J=7$ and $J=8$ by Newton's method}}	
\label{PPR1Table6}											
\end{table}

\subsection{Example 4}

Consider (\ref{eq10}) with $\sigma=-1$,$\gamma=\frac{-1}{2}$ and $\beta=\frac{3}{2}$ along with the boundary conditions,
\begin{equation} \label{eq25}
y(0)=1, ~~y'(1)=y(1).
\end{equation}
This problem is considered in \cite{Implementation1980} to compute various atomic states.

Now, by applying Haar wavelet method (\cite{Haar_solving_HODE2008}), we arrive at the following nonlinear system:
\begin{equation} \label{eq33}
y(t)=1-t+\sum_{i=1}^{2M} a_i (P_{2,i}(t)-t P_{2,i}(1))+x \Big(\sum_{i=1}^{2M} a_i  h(1)\Big)^\frac{2}{3}.
\end{equation}
In this case we compute residual error. Residual error for the differential equation $y''(t)- t^{-1/2} y^{3/2}=0$ is defined as $R(t)=\tilde{y}''(t)- t^{-1/2} \tilde{y}^{3/2}$, where $\tilde{y}$ is approximate Haar wavelet solution. Maximum absolute residual error is defined as $$R^\infty=\max_{i}\left|R(t_i)\right|.$$

\begin{figure}[H]
\begin{center}
\includegraphics[scale=0.3]{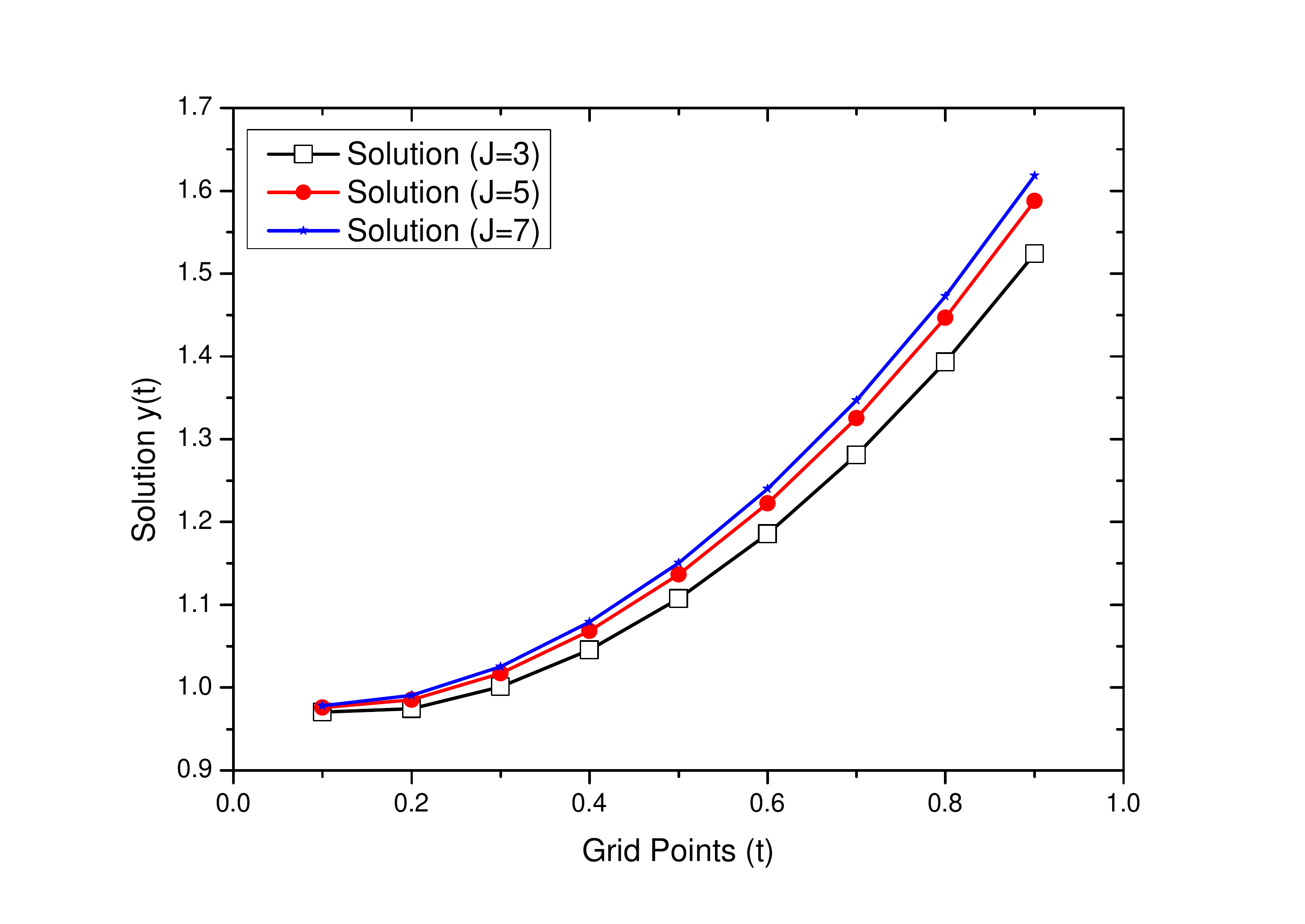}
\caption{Graph of solutions by Newton's method for $J=3$, $J=5$ and $J=7$.}\label{ppr1fig7}
\end{center}
\end{figure}

\begin{figure}[H]
\begin{center}
\includegraphics[scale=0.3]{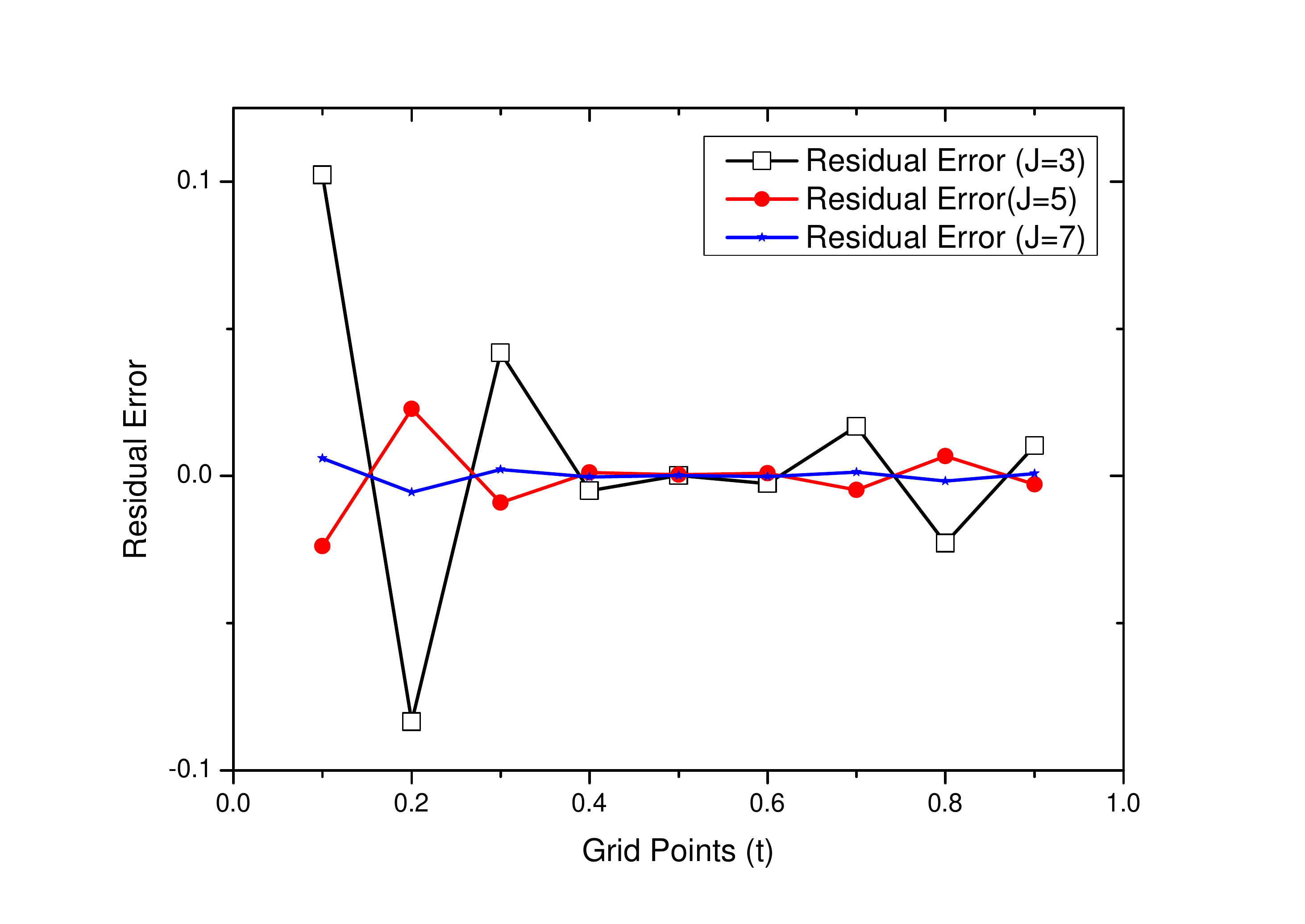}
\caption{Graph of Residual error for Newton's method at $J=3$, $J=5$ and $J=7$.}\label{ppr1fig8}
\end{center}
\end{figure}

To solve nonlinear system obtained by substituting \eqref{eq33} into considered Thomas Fermi equation \eqref{eq10} - \eqref{eq25} with $\sigma=-1$, $\gamma=\frac{-1}{2}$ and $\beta=\frac{3}{2}$, we apply Newton's method and compute the solutions. Taking $[1,1,\cdots,1]$ as initial guess, we get solutions which is depicted in Table \ref{PPR1Table7}. Graph of computed solution and plots for residual errors for $J=3$, $J=5$ and $J=7$ is given in figure \ref{ppr1fig7} and \ref{ppr1fig8}, respectively. To check the continuous dependence on data and stability of the proposed method we take different initial guesses $[1.01,1.01,\cdots,1.01]$ and $[1.1,1.1,\cdots,1.1]$, and observe that final computed solution does not vary significantly. 

As $J$ increases $R^\infty$ decreases and solution is further improved (see Fig. \ref{ppr1fig8}). Since final solution is independent of initial guesses and $R^\infty$ decreases with increase in $J$, our method is a powerful technique which can be used to compute the solution of Emden-Fowler or Thomas Fermi type equations. 

\begin{table}[H]											 
\centering											
\begin{center}											
\resizebox{9cm}{2.5cm}{
\begin{tabular}	{|c | c|  c|  c|}		
\hline
$t$	&	$J=3, N=3$	&	$J=5, N=3$	&	$J=7, N=3$	\\\hline
0.1	&	0.970377	&	0.975761	&	0.978282	\\\hline
0.2	&	0.974685	&	0.985461	&	0.990554	\\\hline
0.3	&	1.00098	&	1.01744	&	1.02527	\\\hline
0.4	&	1.04576	&	1.06831	&	1.07909	\\\hline
0.5	&	1.10753	&	1.13672	&	1.15072	\\\hline
0.6	&	1.18584	&	1.22232	&	1.23985	\\\hline
0.7	&	1.28087	&	1.32541	&	1.34684	\\\hline
0.8	&	1.39328	&	1.44677	&	1.47256	\\\hline
0.9	&	1.5241	&	1.58762	&	1.61829	\\\hline
$R^\infty$	&	0.102233	&	0.0237877	&	0.00601627	\\\hline

\end{tabular}}								
\end{center}
\caption{\small{Solutions for $\sigma=-1$,$\gamma=\frac{-1}{2}$ and $\beta=\frac{3}{2}$ for $J=3$, $J=5$ and $J=7$ by Newton's method.}}	
\label{PPR1Table7}											
\end{table}

\section{Conclusion}\label{Conclusion}

In this paper we have used Haar wavelet collocation method on a class of Emden-Fowler equation which originates in Physics. These equations are also referred as Thomas-Fermi equations. We have proposed two methods: in one approach we first linearise it by Newton's quasilinearization method and then solve it by Haar wavelet collocation method. In the second approach we use Haar wavelet collocation method an arrive at system of non-linear equation which we then solve by Newton-Raphson method. We observed that as value of $J$ is increased beyond certain value computed solution does not vary significantly. We also observed that small perturbation in initial guess does not vary final computed solution at all. So therefore both methods are stable and robust. Higher the resolution we use, our computed solution is closer to the exact solution.

\bibliography{MasterSNarendra}
\bibliographystyle{plain}
\end{document}